\documentclass{amsart}
\usepackage{amssymb}
\usepackage{amsfonts}
\usepackage{amsmath}
\usepackage{latexsym}
\usepackage{hyperref}
\usepackage{mathtools}
\usepackage{tikz, float}
\usepackage[british]{babel}
\usepackage{enumerate}

\usepackage[all]{xy}

\newcommand{\nats}{\mbox{\( \mathbb N \)}}

\newtheorem{theorem}{THEOREM}[section]
\newtheorem{lemma}[theorem]{LEMMA}
\newtheorem{corollary}[theorem]{COROLLARY}

\newtheorem{definition}[theorem]{DEFINITION}

\makeatletter
\let\orgdescriptionlabel\descriptionlabel
\renewcommand*{\descriptionlabel}[1]{%
  \let\orglabel\label
  \let\label\@gobble
  \phantomsection
  \edef\@currentlabel{#1}%
  \let\label\orglabel
  \orgdescriptionlabel{#1}%
}
\makeatother

 \def\set#1{{\{ #1\}}}
 
  \def\nodes{{\sf nodes}}
  
   \def\c#1{{\mathcal #1}}
 
 \def\y{{\sf y}}
 \def\g{{\sf g}}
 \def\r{{\sf r}}
 \def\bb{{\sf b}}% \b mean bold
 \def\w{{\sf w}}
 
 \def\ws{{ winning strategy}}
 
 \def\R#1#2{{\langle   #1^{\conv}\circ#2\rangle}}

 \def\conv{\smile}
\DeclareMathOperator{\comp}{;}

 \def\G{{\bf G}}

\title{A corrected strategy for proving no finite variable axiomatisation exists for RRA}
\author{Rob Egrot} 
\author{Robin Hirsch}
\subjclass[2020]{Primary 	03G15; Secondary 05C90}
\keywords{RRA, Representable relation algebras, Finite variable axiomatisation} 
 
\begin{document}

\begin{abstract} We show that if for all finite $c$ there is a pair of non-isomorphic finite digraphs satisfying some additional conditions, one of which is that they cannot be distinguished in a certain  $c$-colour node colouring game, then there can be no axiomatisation of the class of representable relation algebras in any first-order theory of arbitrary quantifier-depth using only finitely many variables.     This corrects the proposed strategy of Hirsch and Hodkinson, \emph{Relation algebras by games}, North-Holland (2002), Problem 1.  However, even for $c=2$, no pair of non-isomorphic graphs indistinguishable in the game is currently known.
\end{abstract}
\maketitle

\section{Introduction} 
It is known that RRA cannot be axiomatised by any finite theory \cite{Mon64} nor by any equational theory using only finitely many variables \cite[theorem~3.5.6]{Jon88}.  Moreover, any axiomatisation of RRA must involve infinitely many non-canonical equations \cite{HVCan}. To prove that RRA cannot be axiomatised by any $c$-variable   theory would yield the first two of these results as corollaries and would significantly strengthen what is known. 

This note arises from a difficulty with \cite[problem 1, page 625]{HH:book}, which proposed a strategy for approaching the above problem. The proposal, which is mentioned in \cite[p491]{S-A05}, is to find graphs $G$ and $H$ with no homomorphism from $G$ to $H$, but indistinguishable in a certain $c$-colour graph game.   It is claimed that such graphs could be used to prove that there is is no $c$-variable axiomatisation of RRA (for finite $c$). Readers of the book, who failed to solve the problem, cannot be rebuked  since there was a  flaw in the proposed solution. The problem is as follows.

Suppose we can find graphs $G, H$, with no homomorphism from $G$ to $H$, that are indistinguishable in a certain colouring game using $c$ colours (call them \emph{$c$-indistinguishable} for short).  The problem statement uses a `rainbow relation algebra' $\c A_{G, H}$ built from the two graphs.  It can be shown that $\c A_{H, H}$ is representable but $\c A_{G, H}$ is not (see \cite[theorem 16.5]{HH:book}).  It is claimed in the problem statement that, since $G, H$ are $c$-indistinguishable, it follows that $\c A_{H, H}$ and $\c A_{G, H}$ could not be distinguished by any $c$-variable formula (i.e. that $\c A_{G, F}\equiv^c \c A_{H, F}$).  Unfortunately, that implication is false.  The problem is that  there are white atoms $\w_S\in\c A_{G,H}$ for every set of $G$-nodes $S$ of size at most two, which would be represented by binary predicates over graph nodes, but only monadic predicates are used in the graph colouring game.  Thus the proof cannot be completed.

To fix that, the idea is to  let $\c B_{G, H}$ be obtained from $\c A_{G, H}$ by deleting all white atoms $\w_S$ (and deleting any forbidden triple involving these deleted atoms).  This solves one problem, because now it is true that  if $G$ and $H$ are $c$-indistinguishable  then it can be shown  $\c B_{G, F}\equiv^c \c B_{H, F}$.  But it creates another, since $\exists$ really needed atoms $\w_S$ in her \ws\ for the representation game over $\c A_{H, H}$, so we can no longer be sure that $\c B_{H, H}$ is representable.  However, if it happens that every partial homomorphism from $H$ to $H$ of size  two extends to a homomorphism, then  it can be shown that $\c B_{H, H}$ is representable. 

\section{A rainbow construction}\label{S:rainbow}
Let $G, H$ be structures in a signature consisting of only binary predicates, which we refer to as \emph{binary structures}.  For most purposes we can assume $G$ and $H$ to be directed graphs, but it will be convenient later to be able to force graph homomorphisms to preserve non-edges, and so we phrase our results here in terms of binary structures so we can formally handle this without issue.  
 
Given two binary structures $G$ and $ H$, we define an atomic relation algebra $\c B_{G, H}$ by defining its atom structure.    A partial homomorphism is a partial map $h$ from $G$ to $H$ such that if $i\neq i'\in G$ and $(i, i')$ belongs to a binary predicate interpreted in $G$ then $(h(i), h(i'))$ also belongs to that predicate interpreted in $H$.  
The atoms are
\[ \set{1', \bb, \w, \y}\cup\set{ \g_i,: i\in G}\cup\set{ \r_{j,j'}:    j, j'\in H}\]
The non-identity atoms are considered to be black, white, yellow, green or red.  
All atoms are self-converse, except $\r_{j,j'}^\smile = \r_{j',j}$.
Forbidden triples of atoms are Peircean transforms of
\begin{enumerate}
\renewcommand{\theenumi}{\Roman{enumi}}
\item $(1', a, b)$ where $a\neq b$\label{f:1}
\item $(\g_i, \g_{i'}, \g_{i''}), (\g_i, \g_{i'}, \w)$, any $i, i', i''\in G$\label{f:nog}
\item $(\y, \y, \y), (\y, \y, \bb)$\label{f:noy}
\item $(\r_{j_1, j_2}, \r_{j_2', j_3'} , \r_{j_1^*, j_3^*})$, unless $j_1=j_1^*, \; j_2=j_2',\;j_3'=j_3^*$\label{f:red}.
\item \label{f:pim} $(\g_i, \g_{i'}, \r_{j, j'})$ unless $\set{(i, i'), (j, j')}$ is a partial homomorphism.
\item \label{f:ii} $(\g_i, \g_i, \r_{j, j'})$ any  $j,  j'\in H$.

\end{enumerate}
Observe for later that the only combinations of three colours where some but not all triples of atoms with those colours are forbidden are  green-green-red and red-red-red, see \eqref{f:red}, \eqref{f:pim} and \eqref{f:ii}.  

  The relation algebra $\c B_{G, H}$ is the complex algebra of this atom structure. Note that here a triple $(a,b,c)$ being forbidden corresponds to setting $(a\comp b)\cdot c = 0$ in the complex algebra. This is slightly different from the approach taken in \cite{HH:book}, but the differences are entirely superficial.

To build complete representations for these algebras, we will employ a game played by two players, $\forall$ and $\exists$, over a kind of labeled complete digraph. We describe this game now, basing our exposition on the material in \cite[section~11]{HH:book}. An \emph{atomic network} $N$ for an atomic relation algebra $\c A$ consists of a  set of nodes (denoted $\nodes(N)$) and a map (also denoted $N$) from pairs of nodes to atoms of $\c A$, such that
$N(x, x)\leq 1', \; 
N(y, x)= N(x, y)^\conv$ and $(N(x, y), N(y, z), N(x, z))$ is not forbidden, for all $x, y, z\in \nodes(N)$. This property of not containing edge labels forming a forbidden triple $(N(x, y), N(y, z), N(x, z))$ is often referred to as the \emph{consistency} of $N$. For atomic networks $M, N$ we write $M\subseteq N$ if $\nodes(M)\subseteq \nodes(N)$ and for all $x, y\in\nodes(M)$ we have $N(x, y)=M(x, y)$. 

Given an atomic relation algebra $\c A$, the complete representation game for 
$\c A$ has $\omega$ rounds. $\exists$ is trying to build an atomic network for $\c A$, and $\forall$ is trying to force a situation where this is impossible. In a play of this game, let the current atomic network be $N$. Then  $\forall$ picks nodes $x, y\in\nodes(N)$ and atoms $\alpha, \beta$ such that $(\alpha, \beta, N(x, y))$ is not forbidden. In response, $\exists$ is required to extend $N$ to $N'$ such that there is a node $z\in\nodes(N)'$ where $N(x, z)=\alpha,\;N(z, y)=\beta$. The difficulty is that she must label all new edges induced by adding $z$ without causing the resulting network to be inconsistent. In other words, without creating any triangles corresponding to forbidden triples (we also call these \emph{forbidden triangles}).  We can assume that no suitable witness $z$ is already in $N$, else the move is trivial as $\exists$ does not need to add any extra nodes to the network. If $\forall$ has no non-trivial move to make then $\exists$ wins, as the network now reveals the required complete representation.   The game starts with $\forall$ playing a non-identity atom $\alpha$, and $\exists$ creating a two node network $\{x_0,y_0\}$ such that the edge $(x_0,y_0)$ is labeled by $\alpha$. We say $\forall$ wins if in some round he makes a move such that $\exists$ cannot extend the network consistently, and we say $\exists$ wins if she survives $\omega$ rounds, or if $\forall$ cannot make a non-trivial move at some point. We say $\exists$ has a \ws\ if she can play so that her victory is guaranteed. The key result, as proved in \cite[theorem 11.7]{HH:book}, is that an atomic relation algebra with a countable number of atoms $\c A$  is completely representable if and only if $\exists$ has a winning strategy in the complete representation game over $\c A$.

 We return now to the relation algebra $\c B_{G, H}$ defined above.  Given an atomic  network $N$ for $\c B_{G, H}$ and nodes $x, y\in\nodes(N)$ let
\[R_N(x, y)=\set{z\in\nodes(N): N(x, z)\mbox{ is green and } N(y, z)=\y}.\]
Observe that $R_N(x, y)$ depends only on the green and yellow edge labels of $N$.
  A set of nodes of a network where every edge  between distinct nodes has a  red label is called a red clique.  For any $x, y$, by forbidden triples \eqref{f:nog} and \eqref{f:noy}, $R_N(x, y)$ is a red clique.  In a red clique $C$ of size at least two,   by forbidden triple \eqref{f:red}, each node $z\in C$ has a well-defined index $\rho_C(z)\in H$ such that $N(z_1, z_2)=\r_{\rho_C(z_1), \rho_C(z_2)}$, for $z_1\neq z_2\in C$.    So $\rho_C$ is defined on $z\in R_N(x,y)$ by  taking the first  subscript in the label of  $(z, z')$ where $z'\in R_N(x, y)\setminus\set{z}$ is arbitrary. By rule \eqref{f:red}, this subscript does not depend on choice of $z'$.  By consistency of $N$ and \eqref{f:ii}, for each $i\in G$ there can be at most one node $z$ such that $z\in R_N(x, y)$ and $N(x, z)=\g_i$.

 Similarly, if $\theta$ is a complete representation of $\c B_{G,H}$ over base $X$, then for $x, y\in X$, let
 \[ R_\theta(x,y)=\set{z\in X: (x,z)\in\bigcup_{i\in G}\g_i^\theta\wedge (z, y)\in\y^\theta}.\]
  As with networks, if $|R_\theta(x, y)|>1$, then each point $z\in R_\theta(x, y)$  has an index $\rho_{(\theta, x, y)}(z)\in H$. This is defined by noticing that if $z_1\neq z_2 \in R_\theta(x, y)$, then $(z_1,z_2)\in \g_i^\theta\comp \g_{i'}^\theta\;\cap\;\y^\theta;\y^\theta$ for some $i,i'\in G$. As $\theta$ is a complete representation, there is an atom $\alpha$ with $(z_1,z_2)\in \alpha^\theta$  (see \cite[theorem 2.21]{HH:book}) and from the forbidden triple rules we see that $\alpha$ must be $\r_{j,j'}$ for some $j,j'\in H$. We define $\rho_{(\theta, x, y)}(z_1)$ to be $j$, which does not depend on the choice of $z_2\in R_\theta(x, y)\setminus\set{z_1}$.

\begin{theorem} \label{thm:rainbow}
Let $G, H$ be binary structures.
The following are equivalent.
\begin{enumerate}
\item  For all $i\neq i'\in G$ there are $j, j'\in H$ such that $\set{(i, i'), (j, j')}$ is a partial homomorphism, and every partial homomorphism  $\set{(i, j), (i', j')}$ where $i\neq i'$  from a substructure of $G$ into $H$ extends to  a homomorphism $G\rightarrow H$.
\item $\c B_{G, H}$ is completely representable.
\end{enumerate}
\end{theorem}

\begin{proof}

Suppose $\c B_{G, H}$ is  completely representable, say $\theta$ is a complete representation.  Since $\theta$ is complete, for every pair of points $(x, y)\in 1^\theta$ in the base of the representation there is a unique atom $\alpha$ such that $(x, y)\in\alpha^\theta$ (by \cite[theorem 2.21]{HH:book}).  

Let $i_1\neq i_2\in G$. 
   Find points $x, y$ in the base of the representation such that $(x, y)\in\w^\theta$, see  the first part of figure~\ref{fig}.  Since $(\g_{i_t},\y,\w)$ is not forbidden (for $t=1,2$), there are points $z_{1}, z_{2}$ such that $(x, z_{t})\in\g_{i_t}^\theta$ and $(z_{t}, y)\in\y^\theta$, for $t=1,2$.  The unique atom that holds on $(z_{1}, z_{2})$ cannot be the identity by forbidden triple \eqref{f:1}, nor green, white, yellow or black, by forbidden triples \eqref{f:nog}, \eqref{f:noy}, hence it must be red, say $\r_{j_1,j_2}$. We also have $(z_1,z_2)\in \g_{i_1}^\theta\comp \g_{i_2}^\theta$, and so, by forbidden triple \eqref{f:pim}, the map $\set{(i_1,j_1),(i_2,j_2)}$ is a partial homomorphism of size two.  

To show that partial homomorphisms of size two extend,  let  $\set{(i_1, j_1), (i_2, j_2)}$  be a partial homomorphism from $G$ to $H$, where $i_1\neq i_2$.  
Let $z_{i_1}, z_{i_2}$ be distinct points in the base of the representation such that $(z_{i_1}, z_{i_2})\in \r_{j_1, j_2}^\theta$ (see  the second part of figure~\ref{fig}).  Since $\set{(i_1,j_1),(i_2,j_2)}$ is a partial homomorphism it follows from rule \eqref{f:pim} that $(\g_{i_1}, \g_{i_2}, \r_{j_1, j_2})$ is not forbidden, so there is a point $x$ where $(x, z_{i_1})\in\g_{i_1}^\theta$ and $(x, z_{i_2})\in \g_{i_2}^\theta$.  Also, $(\y, \y, \r_{j_1, j_2})$ is not forbidden, so there is a point $y$ where $(y, z_{i_1}), (y, z_{i_2})\in\y^\theta$, and clearly $x\neq y$.  Since the representation is complete, there must be an atom $\alpha$ such that $(x, y)\in\alpha^\theta$ (by \cite[theorem 2.21]{HH:book}), and as $x\neq y$ this atom cannot be $1'$ (by rule \eqref{f:1}).

We have shown that $z_{i_1}\neq z_{i_2}\in R_\theta(x, y)$.  Write $\rho$ for $\rho_{(\theta, x, y)}$, so for $w\in R_\theta(x, y)$, $\rho(w)$ denotes the index of $w$ in $H$, and for $w, w'\in R_\theta(x, y)$ we have $(w, w')\in\r_{\rho(w), \rho(w')}^\theta$.    Since $(z_{i_1}, z_{i_2})\in\r_{j_1, j_2}^\theta$, 
we have $\rho(z_{i_1})=j_1,\;\rho(z_{i_2})=j_2$.

Regardless of which non-identity atom $\alpha$ is, for each node $i$ of $G\setminus\set{i_1,i_2}$, the triple $(\g_{i}, \y, \alpha)$ is not forbidden, from which it follows that $(x,y)\in \g_i\comp \y$, and so there must be a point $z\in R_\theta(x, y)$ where $(x, z)\in g_{i}^\theta$ and $(z, y)\in\y^\theta$.   This point is unique, as if $z'$ is a point with the same properties, then we have $(z,z')\in (\g_i^\theta\comp \g_i^\theta)\cap (\y^\theta\comp \y^\theta)$. As $\theta$ is complete, $(z,z')$ is contained in the interpretation of some atom, and the forbidden triple rules imply that this atom must be the identity. Thus $z=z'$, by definition of the identity in proper relation algebras. The map from $G$ to $H$  that sends $i\in G$ to $\rho(z_i)$ is therefore well defined,
is a homomorphism, by \eqref{f:pim}, and extends $\set{(i_1,j_1), (i_2,j_2)}$ as required.

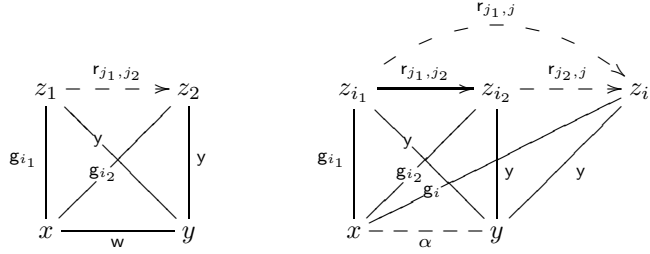
\begin{figure}
\[ 
\xymatrix@!=2.5pc{
z_1\ar@{-}[rd]|-(.3){\y}\ar@{-->}[r]^{\r_{j_1, j_2}}&z_2\ar@{-}[d]^{\y}\\
x\ar@{-}[r]_{\w}\ar@{-}[u]^{\g_{i_1}}\ar@{-}[ur]|-(.4){\g_{i_2}}&y
}\hspace{.5in}
\xymatrix@!=2.5pc{
z_{i_1}\ar@{-}[rd]|-(.3){\y}\ar@{--}@{->}[r]^{\r_{j_1, j_2}}\ar@{-->}@/^2pc/ [rr]^{\r_{j_1,j}}&z_{i_2}\ar@{-}[d]^(.6){\y}\ar@{-->}[r]^{\r_{j_2,j}}&z_i\\
x\ar@{--}[r]_{\alpha}  \ar@{-}[u]^{\g_{i_1}}\ar@{-}[ur]|-(.4){\g_{i_2}}\ar@{-}[rru]|-(.25){\g_{i}}&y\ar@{-}[ru]_{\y}&
}
\]
\caption{\label{fig}From the  representation, a partial homomorphism $\set{(i_1, j_1), (i_2, j_2)}$ exists and extends to $i$.}
\end{figure}

\medskip

Now we check the converse.  Assume the first condition in the theorem. As discussed above, it is sufficient to show she has a \ws\ in the complete representation game for $\c B_{G,H}$.  Assume also that $\forall$ does not play any trivial moves, by which we mean moves where $\exists$ can respond by letting $z$ be a node already in the network. In particular, assume $\forall$ does no play any moves where either $\alpha$ or $\beta$ is $1'$. A consequence of this and $\exists$'s strategy, which we will describe shortly, is that the label $1'$ will only occur in the networks produced during the game as labels for the reflexive edge of a node with itself. 

The basic idea behind $\exists$'s strategy is that she will, as far as possible, use labels which obviously do not interfere with the consistency of the network, and never any labels that are not either white, black or red. The difficult cases turn out to be when she is forced to use a red atom. This occurs only when $z$ is in $R_N'(x,y)$  (or $R_{N'}(y, x)$), and $|R_{N'}(x,y)|>1$ (respectively, $|R_{N'}(y,x)|>1$).

There are three ways to extend the network $N$ to $N'$ so that a new red clique of form $R_{N'}$ is created with $|R_{N'}(x,y)|>1$. The first is if $|R_N(x,y)|=1$ and $\forall$ plays $(x,y,\g_{i'},\y)$ for some appropriate $i'\in G$.  
In this case, $R_{N'}(x, y)=\set{w, z}$ for some $w\in N$, with $N(x, w)=\g_{i}$ for some $i\neq i'\in G$, and by the first part of our assumption there are $j, j'\in H$ such that $h= \set{(i, j), (i', j')}$ is a partial homomorphism.  Here $\exists$ sets $N'(w, z)=\r_{j, j'}$.    
The second way such a red clique of size greater than one can be created is where $w\in N$, \/ $z$ is the new node, and $\forall$'s move is $(x,y,\y,\y)$, with $N(w,x)= \g_i$ and $N(w,y) = \g_{i'}$ for some $i\neq i'\in G$.  Since $\exists$ does not use green or yellow labels, $R_{N'}(w,z)=\set{x, y}$.  Here the map $h$ from green subscripts to red indices is a partial homomorphism of size two, by consistency of the previous network and of $\forall$'s move.  The third way is similar to the second, except $\forall$ plays $(x,y,\g_i,\g_{i'})$, and $N(w,x)=N(w,y)=\y$. Here a partial homomorphism $h$ is defined as in the second case.

In all three cases, the second part of our assumption tells us that $h$ extends to a homomorphism $h^+:G\rightarrow H$.  In later rounds, whenever a new node is added to $R_N(x, y)$, \/ $\exists$ will use $h^+$ to get the index of the new node and hence to label all edges in $R_N(x, y)$ incident with the new node $z$. The fact that  red labels are defined by node indices will ensure that \eqref{f:red} is not violated and the fact that $h^+$ is a homomorphism will ensure that \eqref{f:pim} is not violated. Although in general it turns out that two distinct red cliques can intersect in up to two points, crucially the new point $z$ belongs to at most a single clique in the round when it is added,  so, $\exists$ is never conflicted about which homomorphism to use, as we shall see.

Now we have sketched out $\exists$'s plan, we must check that it works. Since the initial round presents no difficulties, we consider later rounds.   Suppose the current network is $N$ with $|N|\geq 2$, and $\forall$'s move is $(x, y, \alpha, \beta)$, with $\exists$ responding by adding a new node $z$ to create $N'=N\cup\{z\}$.  We have two induction hypotheses: 
\begin{enumerate}
\renewcommand{\labelenumi}{\text{\theenumi}}
\renewcommand{\theenumi}{(\text{H\arabic{enumi}})}
\item\label{eq:induct}  For all $u, v\in N$ where $|R_N(u, v)|\geq 2$ there is a homomorphism $h_{uv} :G\rightarrow H$ such that for all $w\neq w'\in R_N(u, v)$ we have 
\begin{equation*}
 (N(u, w)=\g_i\wedge N(u, w')=\g_{i'})\rightarrow N(w, w')=\r_{h_{uv}(i), h_{uv}(i')}.
 \end{equation*}
\item\label{eq:induct2}  For any $u\neq v\in N$, if $\alpha$ and $\beta$ are green or yellow then there is at most a single node $w\in N$ such that $N(u, w)=\alpha,\; N(w, v)=\beta$.
\end{enumerate}

Given $N$, for each $x', y'\in N$ where $|R_N(x', y')|\geq 2$ let $h_{x'y'}$ be a homomorphism satisfying \ref{eq:induct}.
For each $w\in N\setminus \{x,y\}$ she must assign $N'(w,z)$ in such a way that $N'$ is a consistent network, and the induction hypothesis is maintained.    She proceeds as follows:
\begin{enumerate}[(a)]
\item
 If  $N(w, x)$ and $\alpha$ are not both green, and $N(w, y), \beta$ are not both green, she lets $N'(w, z)=\w$. 
  \item If $N(w, x), \alpha$ are both green but $N(w, y),\beta$ are not both yellow, or if $N(w, y), \beta$ are both green but $N(w, x), \alpha$ are not both yellow, she lets $N'(w, z)=\bb$. 
  \item \label{case:rem} The remaining case is where $N(x, w)=\g_{i},\;\alpha=\g_{i'},\; N(w, y)=\beta=\y$ (or similar with $x, y$ swapped). Note that $i\neq i'$, by the `no trivial moves' assumption.  Here $z,w\in R_{N'}(x,y)$, and she is forced to choose $N'(w, z)=\r_{j,j'}$ for some $j, j'\in H$.  If $h_{xy}$ is already defined  for $N$ she lets $N'(w, z)=\r_{h_{xy}(i), h_{xy}(i')}$, thereby maintaining \ref{eq:induct} for $(x, y)$.  Otherwise, $R_N(x, y)=\set w$ and she may pick any  $j, j'\in H$ such that $\set{(i, j), (i', j')}$ is a partial homomorphism and extend it to a homomorphism $h_{xy}:G\rightarrow H$  (using both parts of  (1)) and again let $R_{N'}(w, z)=\r_{h_{xy}(i), h_{xy}(i')}$, establishing \ref{eq:induct} for  $(x, y)$ in $N'$. Note that if $z \in R_{N'}(x',y')$, then, as $\exists$ never uses green or yellow labels, it's easy to show that $x=x'$ and $y=y'$, so the above strategy is well defined.  
	\end{enumerate}
	
  First we show that $N'$ is a consistent network by checking that the labeling of each triangle $(w, w^*, z)$ for $w\neq w^*\in N$ is not forbidden.  If $\set{w, w^*}=\set{x, y}$ then the triangle is consistent (else $\forall$'s move would be illegal), so without loss of generality we assume that $w\notin \{x,y\}$.  
  Observe that $N'(w, z)$ must be either white, black or red, as, if $N'(w, z)$ is green or yellow it follows that $w\in\set{x, y}$, contradicting our assumption.  If $N'(w, z)=\w$, then the only possibility that the triangle $(w, w^*, z)$ could be forbidden comes from \eqref{f:nog}, but this requires that $N(w^*,z)$ and $N(w,w^*)$ be green, and thus that $w^*\in\{x,y\}$. But then the conditions of (a) would not have been met, so $N'(w,z)$ could not be $\w$ after all.  Similarly, if $N'(w, z)=\bb$ then the possibility of violating \eqref{f:noy} is ruled out by  case (b)  conditions.    
	
	In the remaining case, $N'(w, z)$ is red, and the only forbidden triples involving red atoms are  \eqref{f:red}, \eqref{f:pim} and \eqref{f:ii}.    We assume $N(x, w)=\g_i,\;\alpha=\g_{i'}$, and the case where $x, y$ are swapped follows by symmetry.      A triangle $(w, w^*, z)$ could only violate forbidden triple \eqref{f:red} if all three edges were red, which only happens when $w, w^*, z\in R_{N'}(x, y)$. In this case, by \ref{eq:induct} we have $N'(w,w^*)= \r_{h_{xy}(i), h_{xy}(i^*)}$ for some $i^*\in G$ and some $h_{xy}$, and according to $\exists$'s strategy the other edge labels are $N'(w, z)=\r_{h_{xy}(i), h_{xy}(i')}$ and $N'(z, w^*)=\r_{h_{xy}(i'), h_{xy}(i^*)}$, so \eqref{f:red} is not violated. 
  For forbidden triple \eqref{f:pim},  the only possible  green-green-red triangle incident with $z$ and $w$ is $(x, w, z)$ and the edge labels  $(\g_i, \g_{i'},\r_{h_{xy}(i), h_{xy}(i')})$ do not violate \eqref{f:pim}, since $h_{xy}$ is a homomorphism.     The only triangle  containing $\set{w, z}$ which could violate \eqref{f:ii} is $(w, x, z)$ (only this can be green-green-red), but in this case $i\neq i'$ (else $w$ is a witness to the current move, contrary to the `no trivial moves' assumption) so \eqref{f:ii} is not violated.   Hence $N'$ is a consistent network.
  
It remains to check the induction hypotheses.  \ref{eq:induct2} is clear, since $\exists$ never adds a new node to the network if a suitable witness is already in $N$.  We check
 \ref{eq:induct}.  Suppose $\exists$ is playing according to the strategy we have described, and she adds $z$ to $N$ to obtain $N'$ in response to a move $(x, y, \alpha, \beta)$ by $\forall$.  We say that a pair $(x', y')$ is \emph{safe} if either $x', y' \in N$ and $R_{N'}(x', y')=R_N(x', y'),$ or $|R_{N'}(x', y')|\leq 1$.   If $(x', y')$ is safe then \ref{eq:induct} is true for $(x', y')$ in $N'$ either trivially because the size of $R_{N'}(x', y')$ is less than two, or inductively, since \ref{eq:induct} is assumed true for $N$.

  We check   induction hypothesis \ref{eq:induct} according to whether $\alpha$ and $\beta$ are green, yellow or neither.  
 \begin{enumerate}[(i)]
\item \label{obs:2}
If $\alpha=\g_i,\;\beta=\y$ then $R_{N'}(x, y)=R_N(x, y)\cup\set z$ and  all other red cliques  are safe.
\item If $\alpha=\y,\;\beta=\g_i$ then $R_{N'}(y, x)=R_N(y, x)\cup\set z$ and all other cliques are safe.
\item  If $\alpha=\g_i,\;\beta=\g_{i'}$ and $N(x, w)=N(w, y)=\y$ for some $w\in N$, then $R_{N'}(z, w)=\set{x, y}$, and, by \ref{eq:induct2}, all other red cliques are safe. 
\item If $\alpha=\beta=\y,\; N(x, w)=\g_i, \; N(y, w)=\g_{i'}$ for some $w\in N$,  then $R_{N'}(w, z)=\set{x, y}$,   and, \ref{eq:induct2},  all other red cliques are safe.
\item Else $\set{\alpha, \beta}\not\subseteq\set{\g_i, \g_{i'},\y}$ (for any $i, i'\in G$) and  all red cliques are safe.
\end{enumerate}
For case (i), if $R_{N'}(x, y)=\set{w, z}$ (some $w\in N$) then this red clique has size two in this round, for the first time.    Say $N(x, w)=\g_{i'}$, where $i'\neq i$ else $w$ is already a witness.    By assumption (1) there is a partial homomorphism defined on $\set{i, i'}$ which extends to a homomorphism, and  her strategy chooses such a homomorphism $h_{xy}$ and lets  $N'(w, z)=\r_{h_{xy}(i'), h_{xy}(i)}$,  as required for the induction hypothesis.    If $|R_{N'}(x, y)|\geq 3$ then $|R_N(x, y)|\geq 2$, so inductively there is already a homomorphism $h_{xy}:G\rightarrow H$ determining red labels in $R_N(x, y)$.  In this case, for each $w\in R_{N'}(x,y)$, her strategy defines $N'(w, z)=\r_{h_{xy}(i'), h_{xy}(i)}$, where $N(x, w)=\g_{i'}$, thereby maintaining \ref{eq:induct} for $(x, y)$ in $N'$.     Case (ii) is similar.

 In case (iii)  when $R_{N'}(z, w)=\set{x, y}$, since it follows from the stated conditions that $(\g_i\comp \g_{i'})\cdot(\y\comp \y) \geq N(x, y)$, we know that $N(x, y)=\r_{j, j'}$ for some $j, j'\in H$, and, by \eqref{f:ii}, we must have $i\neq i'$.  By \eqref{f:pim} and the  legality of $\forall$'s move, $\set{(i,j), (i', j')}$ is a partial homomorphism.  By  assumption (1) this extends to a homomorphism $h_{zw}:G\rightarrow H$, as required by the induction hypothesis. Case (iv) is similar, and case (v) is trivial. So the strategy described above is indeed a winning one for $\exists$. 
\end{proof}

\section{Colouring games and pebble games}\label{S:game}
We now define a vertex colouring game played by $\forall$ and $\exists$ over a pair of binary structures using a finite number of colours, which are used to colour \emph{sets} of vertices, rather than individual vertices. This game is studied in more detail for digraphs in \cite{EgrHirgames}, where it is called a \emph{Seurat game} in reference to the pointillist style of painting. A version of the game also appears in \cite{EgrHirNonFin}, where it is defined for sets (which are implicitly treated as complete digraphs). 
Let $0\leq c<\omega,\; 0\leq n\leq \omega$.
Given two binary structures $G, H$ we define the  $c$-colour, $n$-round  colouring game $\G^c_n(G, H)$ to test equivalence of the binary structures using $c$ monadic predicates.  A $G$-interpretation is a map $\set{0, \ldots, c-1}\to\wp(G)$  to subsets of the vertices of $G$, and an $H$-interpretation is a mapping $\set{   0, \ldots, c-1}\rightarrow\wp(H)$. Intuitively, these maps associate vertices in $G$ and $H$ with the different colours.   A \emph{position} in the game consists of a $G$-interpretation and an $H$-interpretation.    

For $n<\omega$, a play of the game is a sequence of $n+1$ positions $((g_0,h_0),\ldots,(g_n,h_n))$, where each pair $(g,h)$ consists of a $G$-interpretation $g$ and an $H$-interpretation $h$. If $n = \omega$ then the game is an infinite sequence of positions $((g_0,h_0),(g_1,h_1),\ldots)$. The maps $g_0$ and  $h_0$ are both defined to be $t\mapsto \emptyset$ for all $t<c$. In other words, all vertices of both graphs start uncoloured.   

If $n = 0$ then neither player does anything and the game ends immediately with $\exists$ being declared the winner. For $n>0$, at the start of  round $r<n$ the current position is $(g_r,h_r)$. Then $\forall$ chooses $t<c$ and a subset of the vertices of $G$ or a subset of the vertices of $H$, \/ $\exists$ responds with a subset of the nodes of the other binary structure. The intuition here is that $\forall$ is colouring some set of vertices in one of the structures, and $\exists$ is responding by colouring a set of vertices of the other structure with the same colour. If $\forall$ reuses a colour that has already been used, then its previous use is first erased from both binary structures. To reflect the new situation,  the position is updated  to $(g_{r+1},  h_{r+1})$ from $(g_r, h_r)$ by changing $g_{r+1}(t)\subseteq G$ and  $h_{r+1}(t) \subseteq H$ according to these choices. 
	
 A \emph{palette}
$\pi$ is a subset of $\set{0, \ldots, c-1}$.  Given a $G$-interpretation $g$, we may interpret $\pi$ by
\begin{equation*}\label{def:pal}\pi^g=\{x\in G: \forall t<c(x\in g(t)\iff t\in \pi) \}.\end{equation*}
Intuitively, $\pi^g$ tells us which vertices of $G$ are coloured according to $g$ with exactly the combination of colours defined by $\pi$.       
Observe that the set of vertices of $G$ is the disjoint union of the sets $\pi^g$, as $\pi$ ranges over palettes.
A position $(g, h)$ is a win for $\forall$ if  either
\begin{enumerate}[(C1)]  
\renewcommand{\labelenumi}{\text{\theenumi}}
\renewcommand{\theenumi}{(\text{C\arabic{enumi}})}
\item\label{C1} there is a palette $\pi$ where $\pi^g$ is empty but $\pi^h$ is not or the other way round,  or 
\item\label{C2} there are palettes $\pi, \pi'$ and a binary predicate $b$ such $(\pi^{g}\times(\pi')^{g})\;\cap\; b^G$ is empty but $(\pi^{h}\times (\pi')^{h})\;\cap \; b^H$ is not, or the other way round.
\end{enumerate}

We say that $\forall$ wins in round $k$ if $(g_k,h_k)$ is the first winning position for him.  For $n<\omega$, if $\forall$ does not win in any round $i\leq n$, then $\exists$ is the winner. If $n=\omega$, then $\exists$ wins if $(g_k,h_k)$ is not a win for $\forall$ for all $k<\omega$.

\medskip
In addition to the game described above, we will use the following minor variation of the classic Ehrenfeucht-Fra\"iss\'e game. Given two relation algebras $\c A, \c B$ we define the $c$-pebble, $n$-round equivalence game $\Gamma^c_n(\c A, \alpha_0, \c B,\beta_0)$.  The pair $(\alpha_0, \beta_0)$ defines the starting position of the game, and consists of two partial maps $\alpha_0:\set{0, \ldots, c-1} \rightarrow \c A$ and $\beta_0:\set{0, \ldots, c-1}\rightarrow \c B$. We require that $\alpha_0$ and $\beta_0$ have the same domains (i.e. that they are defined for the same elements). The intuition is that $\alpha_i(t)$ and $\beta_i(t)$ denote the position of the $t$-coloured pebbles at the start of round $i$, and if $t$ is not their domain then that pebble pair has not been placed yet. If the initial position is defined by maps with empty domains, then we may refer to the game just as $\Gamma^c_n(\c A, \c B)$ for brevity.

If $n=0$, then the game is entirely determined by the starting position, and neither player does anything. For $n\geq 1$, in each round $k<n$, if the position is $(\alpha_k, \beta_k)$, \/  $\forall$ picks $t<c$ and an element of $\c A$ or of $\c B$, and then $\exists$ picks an element of the other algebra.  At the end of the round, the position is updated by changing $\alpha_k(t)\in\c A,\;\beta_k(t)\in\c B$ according to these 
choices, but leaving other values unchanged. 
This defines $\alpha_{k+1}$ and $\beta_{k+1}$ ready for the start of the next round. Let $\c A_\alpha$ and $\c B_\beta$ 
denote the subalgebras of $\c A$ and $\c B$ generated by the images of $\alpha$ and $\beta$, respectively.  
 At the start of round $k$, consider the binary relation $\alpha_k^{\conv}\circ\beta_k=\set{(\alpha_k(t), \beta_k(t)): t\in dom(\alpha_k)}$ (here $\circ$ and ${}^\conv$ denote, respectively, ordinary composition and conversion of relations). 

We aim to use $\alpha_k^{\conv}\circ\beta_k$ to define a function $\R{\alpha_k}{\beta_k}: \c A_{\alpha_k}\to \c B_{\beta_k}$. To do this we note that elements of $\c A_{\alpha_k}$ correspond to terms constructed from relation algebra constants and elements of $\alpha_k$ using relation algebra operations. Similar applies to elements of $\c B_{\beta_k}$. We want to define $\R{\alpha_k}{\beta_k}$ so that it takes such a term in $\c A_{\alpha_k}$ to a term in $\c B_{\beta_k}$ by fixing relation algebra constants, and swapping elements of $\alpha_k$ with their partners in $\beta_k$. 

We say $\forall$ wins the game in round $k$ if   $\R{\alpha_k}{\beta_k}$ is not an isomorphism, or if it fails to be a well defined function at all, and the maps have been isomorphisms in 
all previous rounds. On the other hand, $\exists$ wins if the maps we have described are isomorphisms for all $k\leq n$. The value of these modified Ehrenfeucht-Fra\"iss\'e games is given by the following definition and lemma.

	\begin{definition}
	Let $\c A$ and $\c B$ be relation algebras, and let $\alpha$ and $\beta$ be partial maps from $\{0,\ldots,c-1\}$ to $\c A$ and $\c B$, respectively, and suppose also that $\alpha$ and $\beta$ have the same domains. We say 
	\[(\c A, \alpha) \equiv^c_n (\c B, \beta)\]
	if whenever $\phi$ is a first-order formula in the language of relation algebras with the additional restrictions that the quantifier depth in $\phi$ be at most $n$, that $\phi$ involves only  variables from the set $\{x_0,\ldots,x_{n-1}\}$, and that the free variables of $\phi$ are all indexed by values from the domain of $\alpha$ and $\beta$, we have 
	\[\c A, \alpha \models \phi \iff \c B, \beta \models \phi.\]
	Here, for example, $\c A, \alpha \models \phi$ means that $\c A \models \phi$ if all variables $x_i$ occurring free in $\phi$ are assigned to $\alpha(i)$ in $\c A$. When $\alpha$ and $\beta$ are empty we just write $\c A \equiv^c_n \c B$.
	\end{definition}
		
\begin{lemma}\label{L:EF}
Let $1\leq c < \omega$ and $n \leq \omega$, let $\c A$ and $\c B$ be relation algebras, and let $\alpha_0$ and $\beta_0$ be partial maps from $\{0,\ldots,c-1\}$ to $\c A$ and $\c B$, respectively, with the same domains. If $\exists$ has a winning strategy in $\Gamma^c_n(\c A,\alpha_0, \c B, \beta_0)$  then we must have $(\c A, \alpha) \equiv^c_n (\c B, \beta)$. 
\end{lemma}
\begin{proof}
This is half the well known result for relational signatures (see e.g. \cite[theorem 6.10]{Im99}). Having functions in the signature blocks the proof of the converse. We induct on $n$. For the base case, $\exists$ has a winning strategy in $\Gamma^c_0(\c A,\alpha_0, \c B, \beta_0)$ if and only if the induced map $\R{\alpha_0}{\beta_0}$ is an isomorphism, if and only if $(\c A, \alpha_0), (\c B, \beta_0)$ agree on all equations using appropriate variables, if and only if $(\c A, \alpha_0) \equiv^c_0 (\c B, \beta_0)$.

For the inductive step, suppose $\exists$ has a winning strategy in $\Gamma^c_{n+1}(\c A, \alpha_0, \c B, \beta_0)$.   Let $\phi=\exists x_i \psi$ be a formula of quantifier depth at most $n+1$, where the variables occurring free in $\psi$ are either indexed by values for which $\alpha_0$ and $\beta_0$ are defined, or are $x_i$.    If $\c A, \alpha_0\models \exists x_i\psi$ then there is an $x_i$-variant $\alpha_1$ of $\alpha_0$ such that $\c A, \alpha_1\models \psi$.  If $\forall$ plays $\alpha_1(x_i)$ in the game, then since $\exists$ has a winning strategy there is an $x_i$-variant $\beta_1$ of $\beta$ where $\exists$ has a \ws\ in $\Gamma^c_n(\c A, \alpha_1, \c B, \beta_1)$.    Inductively, $\c B, \beta_1\models\psi$, hence $\c B, \beta_0\models\exists x_i \psi$.  Since the argument is symmetric, it follows that $(\c A, \alpha_0)$ agrees with $(\c B, \beta_0)$ on all $c$-variable formulas  $\exists x_i\psi$ where $\psi$ has quantifier depth at most $n$, hence they agree on all $c$-variable formulas of quantifier depth at most $n+1$.      By induction, the lemma holds for all finite $n$.  For the case $n=\omega$ then a \ws\ for $\exists$ in $\Gamma^c_\omega(\c A, \alpha_0, \c B, \beta_0)$ entails a \ws\ in all finite length games, so $(\c A, \alpha_0)\equiv^c_n(\c B, \beta_0)$ for all finite $n$.  Hence $(\c A, \alpha_0)\equiv^c_\omega(\c B, \beta_0)$, as required.
\end{proof}

\section{A corrected strategy and why it looks difficult}\label{S:correction}
It turns out that if $\exists$ has a winning strategy in the infinite game $\G^{c+3}_\omega(G,H)$, then she also has a winning strategy in $\Gamma^c_n(\c B_{G,F}, \c B_{H,F})$ for all $n<\omega$, and this can be used to correct the claims of \cite[problem 1, page 625]{HH:book}. We will prove this soon, but first we will need the following lemmas.

\begin{lemma}\label{L:size}
Let $G$ and $H$ be binary structures and let $c\geq 2$. Then, if $\exists$ is playing $\G^c_\omega(G,H)$ according to a winning strategy, whenever $\forall$ colours a set in one of the binary structures, $\exists$ must respond by colouring a set of nodes of the other binary structure with the same cardinality 
\end{lemma}
\begin{proof} If the set of nodes with a certain colour is bigger in one structure than the other, then $\forall$ may use a second colour to colour all but one node in the larger set (and $\exists$ must colour a proper subset of the smaller set to avoid losing straight away), and he may repeat by re-using his first colour to colour all but one node of the larger set, and so on, until he colours a non-empty set of nodes in the first graph but $\exists$ has only the empty set to choose in the other graph, so $\forall$ wins.  
See \cite[Proposition 2.3]{EgrHirgames} for the details.
\end{proof}

Now consider the relation algebra equivalence game $\Gamma^c_\omega(\c B_{G, F}, \c B_{H, F})$.
We suppose that $\exists$ maintains a private corresponding play of a colouring game over $(G,H)$. Specifically, we suppose she is playing according to a winning strategy in the infinite game $\G^{c+3}_\omega(G, H)$, which has three extra colours.  
So, to recap, if, for example, $\forall$ picks an element $x\in\c B_{G, F}$, she picks the element $y\in\c B_{H, F}$ with the identical non-green part and green part determined by $H_y$, where $H_y$ is the response to $G_x$ in the parallel play of $\G^{c+3}_\omega(G, H)$.    These moves in the play of $\G^{c+3}_\omega(G, H)$ are determined by the play of $\Gamma^c_\omega(\c B_{G, F}, \c B_{H, F})$ and only involve the first $c$ colours.  Our assumption is that  at each position $(g, h)$ occurring in the play of  $\G^{c+3}_\omega(G, H)$,\/   $\exists$ has a \ws\ in the game proceeding from $(g, h)$, even if $\forall$ decides to use the three additional colours. In other words, $\exists$ cannot make any move in the parallel game resulting in a position from which $\forall$ could force a win.

Suppose $(\alpha, \beta)$  is a position in $\Gamma^c_\omega(\c B_{G, F}, \c B_{H, F})$, played as described above, and let $(g,h)$ be the corresponding position in $\G^{c+3}_\omega(G, H)$. We will need to interpret terms in the language of relation algebras with finite variable set $\{x_0,\ldots,x_{c-1}\}$ in the algebras $\c B_{G, F}$ and $\c B_{H, F}$.    For any variable $x_i$ where $i\in dom(\alpha)$, we interpret $x_i$ in $\c B_{G, F}$ by defining $x_i^\alpha=\alpha(i)\in\c B_{G, F}$. If $b$ is a relation algebra constant, we define $b^\alpha$ to be the interpretation of $b$ in $\c B_{G, F}$.  Thus, any relation algebra term  $t$ involving only variables with indices in $dom(\alpha)$ and relation algebra constants has an obvious interpretation $t^\alpha\in\c B_{G, H}$, and similarly $t^\beta\in\c B_{H, F}$.  For any such term $t$ we define $\gamma(t)=\set{x\in G:\g_x\leq t^\alpha}$ and $\eta(t)=\set{y\in H:\g_y\leq t^\beta}$.
\begin{lemma}\label{L:terms}
Let $(\alpha,\beta)$ be a position arrived at in a game $\Gamma^c_n(\c B_{G,F}, \c B_{H,F})$ during which $\exists$ plays using a winning strategy in a parallel game $\G^{c+3}_\omega(G,H)$, as described above. Let $t$ be a term involving only variables indexed by values from $dom(\alpha)=dom(\beta)$. Then:
\begin{enumerate}
\item The sets of non-green atoms below $t^\alpha$ and $t^\beta$ are identical. \label{e:ng}
\item In the parallel game $\G^{c+3}_\omega(G, H)$, if $\forall$ were to use a colour not previously used to colour $\gamma(t)$, then $\exists$   would have to respond by colouring $\eta(t)$, otherwise $\forall$ could force a win, and similar with  $\gamma(t)$ and $\eta(t)$ switched. \label{e:win}
\item $|\gamma(t)|=|\eta(t)|$. \label{e:card}
\end{enumerate}
\end{lemma}
\begin{proof}
For convenience we match the colours $\{0,\ldots,c-1\}$ in $\Gamma^c_n(\c B_{G,F}, \c B_{H,F})$ and $\G^{c+3}_\omega(G,H)$ in the obvious way, and we refer to the colours $\{c,c+1,c+2\}$ used in $\G^{c+3}_\omega(G,H)$ as \emph{additional colours}, or words to that effect.

We will use induction on $t$ to prove \eqref{e:ng} and \eqref{e:win}, and we note that \eqref{e:card} follows from \eqref{e:win}, because if $|\gamma(t)|\neq|\eta(t)|$, then by colouring $\gamma(t)$ with one of the extra colours, $\forall$ could force $\exists$ to colour a set with a different size, and thus force a win in $\G^{c+3}_\omega(G,H)$ (see Lemma \ref{L:size}). In the base case, \eqref{e:ng} is automatic. For \eqref{e:win}, suppose $t = x_i$, and $t^\alpha = \alpha(i)$ for some $i\in dom(\alpha)$. So $\gamma(t) =\gamma(x_i)=\{x\in G:\g_x\leq \alpha(i)\}$ is already coloured by the $i$th colour. Moreover, $\eta(t)=\{y\in H: \g_y\leq \beta(i)\}$ must also be coloured by this colour, as according to $\exists$'s strategy $\beta(i)$ is defined to make this true. If $\forall$ uses one of the additional colours to colour $\gamma(t)$, then $\exists$ must colour all of $\eta(t)$ in response, otherwise there will be a palette mismatch between $G$ and $H$. The cases where $t$ is one of the relation algebra constants are also easy. For the inductive step, we proceed as follows (assuming the result for terms $s,s_1,s_2$):

$t = -s$: As the non-green parts of the two relation algebras are identical, \eqref{e:ng} holds for $t$. For \eqref{e:win}, if $\forall$ uses an additional colour to colour $\gamma(-s)\subseteq G$ then $\exists$ must colour a set $Y\subseteq H$ with the same  colour.  If he goes on to colour $\gamma(s)\subseteq G$ with a second additional colour, then by inductive assumption $\exists$ will colour $\eta(s)\subseteq H$ with that colour.    We know that $\gamma(s), \gamma(-s)$ are disjoint and cover $H$.  Since the position is a winning position it must be that  $Y, \eta(s)$ are disjoint and cover $G$, hence $Y$ is the complement in $H$ of $\eta(s)$, so $Y=\eta(-s)$, as required.
 
$t = s_1\cdot s_2$: For \eqref{e:ng}, if $a$ is a non-green atom, then, appealing to the inductive hypothesis, we have
\[a\leq (s_1\cdot s_2)^\alpha \iff a\leq s_1^\alpha \wedge a\leq s_2^\alpha \iff a\leq s_1^\beta \wedge a\leq s_2^\beta \iff a\leq (s_1\cdot s_2)^\beta.\] 
For \eqref{e:win}, if $\forall$ uses an additional colour to colour $\gamma(s_1\cdot s_1)\subseteq G$ then $\exists$ must colour some set $Y\subseteq H$ with the same colour.  If he went on in the next two rounds to use the other two additional colours to colour $\gamma(s_1)$ and then $\gamma(s_2)$, then inductively we know that $\exists$ colours $\eta(s_1), \eta(s_2)\subseteq H$.  Since the position at the end of this is not a win for $\forall$ we must have $Y=\eta(s_1)\cap\eta(s_2) = \eta(s_1\cdot s_2)$, as required.

$t = s^\conv$: For any non-green atom $a$ we know $a^\conv$ is also a non-green atom and so
\[a\leq (s^\conv)^\alpha\iff a^\conv \leq (s^\alpha)^\conv \iff a^\conv\leq s^\alpha \iff a^\conv\leq s^\beta\iff a\leq (s^\conv)^\beta,\] proving \eqref{e:ng}.  Part \eqref{e:win} is easy since all green elements are self-converse, so $\gamma(s^\conv)=\gamma(s)$ and $\eta(s^\conv)=\eta(s)$.

$t = s_1\comp s_2$: We show first that \eqref{e:ng} holds. If $a\leq (s_1\comp s_2)^\alpha$ is a non-green atom then there are atoms  $b_1\leq s_1^\alpha,\;b_2\leq s_2^\alpha$ where $a\leq b_1\comp b_2$.     Not all triples of atoms of the colours of $b_1, b_2, a$ are forbidden, as we are assuming $a\leq b_1\comp b_2$.   If no triple of atoms of the colours of $b_1, b_2, a$ is forbidden, then $a\leq (s_1\comp s_2)^\beta$, because, by the inductive hypothesis with \eqref{e:ng} and \eqref{e:card}, the non-green atoms below $s_1^\alpha$ and $s_1^\beta$ are the same, and likewise for $s_2$, and $|\gamma(s_1)|=|\eta(s_1)|$, and likewise for $s_2$.

This leaves the case where some but not all triples of atoms of the colours of $b_1, b_2, a$ are forbidden. Recall that the only such colour combinations are red-red-red and green-green-red.    If $b_1, b_2, a$ are all red then, by the inductive hypothesis on \eqref{e:ng}, we have $b_1\leq s_1^\beta,\;b_2\leq s_2^\beta$, so $a\leq b_1;b_2\leq  (s_1;s_2)^\beta$.   This leaves the case where $a=\r_{j_1,j_2}$ is a red atom, and $b_1=\g_{x_1}, b_2=\g_{x_2}$ are both green.    Since $(\g_{x_1},\g_{x_2},\r_{j_1,j_2})$ is not forbidden,  we know that $x_1\neq x_2$ and $\set{(x_1, j_1), (x_2, j_2)}$ is a partial homomorphism.  Since $(g, h)$ is a position from which $\exists$ has a winning strategy in $\G^{c+3}_\omega(G, H)$, after a short but technical argument to come, we shall see that there must be $y_1\in\eta(s_1),\; y_2\in\eta(s_2)$ such that $\set{(y_1, j_1), (y_2, j_2)}$ is a partial homomorphism.

To this end, suppose for contradiction that there are no such $y_1, y_2$.  
In the parallel game $\G^{c+3}_\omega(G, H)$, if  $\forall$ were to colour $\gamma(s_1)$ with one of the extra colours, and then to use another to colour $\gamma(s_2)$, then, by the inductive hypothesis on \eqref{e:win}, $\exists$ would have to respond by colouring $\eta(s_1)$ and $\eta(s_2)$ respectively. Suppose that $\forall$ then uses the final additional colour to colour $\{x_1, x_2\}$. Then $\exists$ must colour some $\{y_1,y_2\}\subseteq H$. Since $x_1\in \gamma(s_1)$ and $x_2\in \gamma(s_2)$,  we can suppose without loss of generality $y_1\in \eta(s_1)$ and $y_2\in \eta(s_2)$. Suppose that, in a further move, $\forall$ uses one of the original colours (i.e. $\{0,\ldots,c-1\}$)  to colour $\set{x_1}$, and without loss of generality we can assume that $\exists$ responds by colouring $\set{y_1}$ (if $\exists$ can colour $\set{y_2}$ then we must have $x_1,x_2\in\gamma(s_1)\cap \gamma(s_2)$ and $y_1,y_2\in\eta(s_1)\cap \eta(s_2)$, in which case we can just switch the labels of $y_1$ and $y_2$).  But then  $\set{(y_1, j_1), (y_2, j_2)}$ is not a partial homomorphism but $\set{(x_1, j_1),(x_2,j_2)}$ is a partial homomorphism, and it follows that $\set{(x_1,y_1),(x_2,y_2)}$ is not a partial isomorphism, indicating that \ref{C2} holds. Thus $\forall$ wins, contrary to our assumption that $\exists$ is following a winning strategy.
    The implication $a\leq (s_1;s_2)^\beta\Rightarrow a\leq(s_1;s_2)^\alpha$, for non-green atoms $a$, is proved similarly. This proves \eqref{e:ng}.

For \eqref{e:win}, the cases where either $s_1^\alpha$ or $s_2^\alpha$ is zero or the identity are trivial, so assume not.    
First we suppose that $s_1^\alpha$ is either a single non-green atom, or pure green (i.e. above only green atoms), and the same for $s_2^\alpha$. By our induction hypothesis, $s_1^\beta$ is the same non-green atom in the former case and green in the latter case, and similar for $s_2^\beta$.  Suppose the colour of $s_1^\alpha$ is $c_1$, and that the colour of $s_2^\alpha$ is $c_2$. 
  If all triples of atoms of colours $c_1$-$c_2$-green are forbidden (i.e. $c_1$ is green, $c_2$ is green or white, or the other way round) then $\gamma(s_1\comp s_2)=\eta(s_1\comp s_2)=\emptyset$.    If no triple of atoms of colours  $c_1$-$c_2$-green 
  are forbidden then 
  $\gamma(s_1;s_2)=G$
   and $\eta(s_1;s_2)=H$.    In both these cases, if $\forall$ colours $\gamma(s_1\comp s_2)$, then $\exists$ must obviously respond by colouring $\eta(s_1\comp s_2)$. 
  The only colour combinations where some but not all triples of atoms of those colours are forbidden are red-red-red and green-green-red.   
   So the remaining cases are where $c_1$ is green, $c_2$ is red, or the other way round.   Without loss of generality, suppose $s_1^\alpha$ is  green and $s_2^\alpha=\r_{j, j'}$ is a red atom.    Observe that
  \[\gamma(s_1\comp s_2)=\set{x'\in G: \exists x\in\gamma(s_1)\text{ s.t. } \set{(x,j), (x', j')}\mbox{ is a partial homomorphism}},\] and 
	\[\eta(s_1\comp s_2)=\set{y'\in H:\exists y\in\eta(s_1)\text{ s.t. } \set{(y,j), (y', j')}\mbox{ is a partial homomorphism}},\] using a Peircean transformation with \eqref{f:pim}. If either of these is empty then the result is trivial, so we assume not.

  If $\forall$ makes an additional move by colouring $\gamma(s_1\comp s_2)\subseteq G$ with $k$ where $c\leq k<c+3$, then $\exists$ responds by colouring some set $Y\subseteq H$ with $k$.    We have to prove that $Y= \eta(s_1\comp s_2)$. 
  
   If he made another additional move by colouring $\gamma(s_1)$ with another additional colour $k'$ then, by our induction hypothesis, $\exists$ would colour $\eta(s_1)\subseteq H$ with that colour. 	
For any $y'\in\eta(s_1\comp s_2 )$, \/  $\forall$ could use the third additional colour $k''$ to colour $\{y'\}$, and in response, $\exists$ would have to colour some $\{x'\}\subseteq G$.  He could continue by picking $y\in \eta(s_1)$ such that $\set{(y, j), (y', j')}$ is a partial homomorphism and colouring $\set{y}$ with one of the original colours (say colour $0$) and $\forall$ would have to colour $\set{x}$ with $0$, where $x\in\gamma(s_1)$ as $\gamma(s_1), \eta(s_1)$ are coloured with $k'$.  Since $\forall$ does not win at this point, we know that $\set{(x, y), (x', y')}$ is a partial isomorphism, hence $\set{(x, j), (x', j')}$ is a partial homomorphism and $x'\in\gamma(s_1\comp s_2)$.  Since $x'$ is coloured with $k$ and $\forall$ does not win, $y'$ must also be coloured $k$, and so $y'\in Y$.

Conversely, suppose $y'\in Y$ (and so is coloured with $k$).    If $\forall$ uses the third additional colour $k''$ for $\set{y'}$ then $\exists$ must respond with $\set{x'}\subseteq \gamma(s_1\comp s_2)$. \/ $\forall$ picks $x\in\gamma(s_1)$ such that $\set{(x, j), (x', j')}$ is a partial homomorphism and colours $\set{x}$ with $0$, then $\exists$ responds by colouring $\set{y}$ with $0$, where $y\in\eta(s_1)$ and  $\set{(y, j), (y', j')}$ is a partial homomorphism.  Hence $y'\in\eta(s_1\comp s_2)$.  Thus $Y=\eta(s_1\comp s_2)$.	
  
More generally, $s_1^\alpha$ and $s_2^\alpha$ are sums of non-green atoms and a single green element, and \eqref{e:win} holds for $(s_1\comp s_2)^\alpha$ since, as $\comp$ is additive, $s_1\comp s_2$ is the sum of simple terms where we have just shown it holds (the inductive case where $t=s_1+s_2$ being covered by the $-$ and $\cdot$ cases).

\end{proof}

 \begin{corollary}\label{C:terms}
If $\exists$ has a winning strategy in $\G^{c+3}_\omega(G,H)$, then she also has a winning strategy in $\Gamma^c_n(\c B_{G,F}, \c B_{H,F})$ for all $n\in\nats$.
\end{corollary}
\begin{proof}
We assume that $\exists$ plays $\Gamma^c_n(\c B_{G,F}, \c B_{H,F})$ using the strategy used in lemma \ref{L:terms}. Suppose $(\alpha,\beta)$ is a position reached during play,  and let $t$ be a term such that $t^\alpha\neq 0$.  Then there is an atom below $t^\alpha$, and so by lemma \ref{L:terms} there must be an atom below $t^\beta$. By this and symmetry we have $t^\alpha =0 \iff t^\beta = 0$, and it follows that $\forall$ does not win at $(\alpha, \beta)$, as the map $\R{\alpha}{\beta}:(\c B_{G,F} )_{\alpha}\to(\c B_{H,F})_{\beta}$, which sends $t^\alpha$ to $t^\beta$, is thus an isomorphism.
\end{proof}

The following theorem, based on our modified rainbow algebras, corrects the argument presented in \cite[problem 1]{HH:book}.

\begin{theorem}\label{T:correction}
Suppose $G, H$ are finite binary structures such that
\begin{enumerate}[1.] 
\item \label{en:strat} $\exists$ has a winning strategy in $\G^{c+3}_\omega(G, H)$, 
\item every partial homomorphism of $H$ of size two extends to a full homomorphism  from $H$ into itself, and 
\item \label{en:no h}  there are $i\neq i'\in G$ and $j,j'\in H$ such that $\{(i,j),(i',j')\}$ is a partial homomorphism that does not extend to a homomorphism $G\to H$.

\end{enumerate}
Then RRA cannot be defined by any $c$-variable first-order theory.
\end{theorem}

\begin{proof}
It follows from lemma \ref{L:EF} and corollary \ref{C:terms} that $\c B_{G, H}$ and $\c B_{H, H}$ would agree about all first-order formulas with at most $c$ variables, but  by theorem~\ref{thm:rainbow} the latter relation algebra would be completely representable while the former would not.  Since both algebras are finite, all representations are complete, and it would follow immediately that
  $\c B_{H, H}  \in RRA,\; \c B_{G, H}\not\in RRA$.

\end{proof}

Note that, in the cases of interest, if there is no homomorphism from $G$ to $H$, then \eqref{en:no h} is satisfied. This is because the negation of \eqref{en:no h} is that for every $i\neq i'\in G$ and $j,j'\in H$, if $\{(i,j),(i',j')\}$ is a partial homomorphism then it extends to a homomorphism $G\to H$.  With no homomorphism from $G\to H$, the only way this can happen  is if every $i\neq i'\in G$ and $j,j'\in H$  it happens that $\{(i,j),(i',j')\}$ is a not partial homomorphism. This can only happen if $G$ is complete and $H$ is edgeless, which is excluded by \eqref{en:strat}.

\begin{corollary}\label{C:digraph}
Let $G, H$ be finite digraphs such that 
\begin{enumerate}[1.]
\item $G, H$ cannot be distinguished in a modified infinitely long $(c+3)$-colour game where $\forall$ can also win at position $(g, h)$  if there are two palettes $\pi, \pi'$ and every pair from $\pi^g\times(\pi')^g$ is an edge but not every pair from $\pi^h\times (\pi')^h$ is an edge, or the other way round,  
\item every partial embedding of $H$ into itself of size two extends to an automorphism of $H$, and
\item there is no embedding of $G$ into $H$.
\end{enumerate}
Then RRA cannot be defined by any $c$-variable first-order theory.
\end{corollary}
\begin{proof}
Two digraphs $G, H$ may be considered as binary structures with three predicates, one for edges, another for non-edges and a third for `non-equality' given by $\{(u,v):u\neq v\}$. The modified game for digraphs understood as binary structures with a single edge relation (the standard setting) is equivalent to the original game played over the same digraphs understood as binary structures with `edge', `non-edge' and `non-equality' relations. To see this, note that for palettes $\pi_1$ and $\pi_2$ the non-edge relation holding between $\pi_1^g$ and $\pi_2^g$ and not between $\pi_1^h$ and $\pi_2^h$ results in a win for $\forall$ according to the modified rules, even if we only consider the `edge' relation. Moreover, when $\pi_1\neq \pi_2$ the `non-equality' relation holds between $\pi_1^g$ and $\pi_2^g$ whenever these are both non-empty, as interpretations of distinct palettes are disjoint. Similar holds for $\pi_1^h$ and $\pi_2^h$, and so accommodating `non-equality' in the standard digraph setting does not require any modification to the game rules.   

Now, homomorphisms in the `three relation' setting clearly correspond to embeddings in the standard digraph setting, and so condition 2 of theorem \ref{T:correction} translates into condition 2 here. Moreover, if there exists an embedding $G\to H$, then, by assumption of condition 2, every partial embedding $G\to H$ of size two must extend to a full embedding of $G$ into $H$, so condition 3 here covers both conditions 3.a and 3.b from theorem \ref{T:correction}.  
\end{proof}

The result of corollary~\ref{C:digraph} and the graph result obtained from theorem \ref{T:correction} by writing `digraph' for `binary structure' seem to be incomparable in strength, though again the uncertainty around condition 1 prevents us from being sure. We get some indication of this by examining the conditions on the graph $H$ in the two results. Respectively, these are:
\begin{enumerate}[(1)]
\item Every partial homomorphism of size two of $H$ into itself extends to a full homomorphism.
\item Every partial embedding of $H$ into itself of size two extends to an automorphism.
\end{enumerate}
Observe that the cyclic graph $C_4$ satisfies (2) but not (1), and we can construct a graph satisfying  (1) but not (2) as follows (we work with undirected graphs here for simplicity). Let $W_2$ and $W_3$ be walk (AKA chain) graphs, and let $R$ be the graph with a single reflexive vertex $v$. Define $H$ by taking the disjoint union of $W_2$, $W_3$ and $R$, and adding an edge $(w,v)$ for each $w\in W_2\cup W_3$. Then $H$ satisfies (1), as every partial homomorphism can be extended by sending every other vertex to $v$. On the other hand, $H$ does not satisfy (2) as, if we suppose the vertices of $W_2$ and $W_3$ are $\{u_0,u_1\}$ and $\{w_0,w_1,w_2\}$ respectively, the partial embedding $\{(w_0,u_0),(w_1,u_1)\}$ cannot be extended to an automorphism.  

Looking at the conditions on $H$ may well be beside the point however, as, as mentioned previously, we do not know whether non-isomorphic graphs $G, H$ exist such that $\exists$ can win even $\G^{2}_\omega(G, H)$.  The power of the game $\G$ to distinguish between graphs is investigated in more detail in \cite{EgrHirgames}. In particular, it is shown there in section 8 that finding non-isomorphic graphs indistinguishable in the 3-colour game would disprove the reconstruction conjecture for graphs (or one of its variations, in the case of directed graphs). This indicates that, at the very least, such graphs will likely be difficult to find.

\bibliographystyle{abbrv}

%\bibliography{../../robin}

 \def\www{/\allowbreak}

\end{document}